# ON THE RELATIVE QUANTITIES OCCURRING WITHIN PHYSICAL DATA SETS

Alex Ely Kossovsky

ABSTRACT

A statistical measure is given expressing relative occurrences of quantities within a given data set. Application of this measure on several real life physical data sets and some abstract distributions are shown to yield quite consistent results. These empirical results also correspond almost exactly to the theoretical converging limit of such a measure mathematically constructed for k over x distribution defined over an infinite range.

## [1] Number System Invariance Principle

That the first significant digits obey the analytical expression Probability[1st digit d] = $LOG_{10}(1 + 1/d)$ has long been accepted as an important probabilistic law prevalent in numerous real life data sets. BL is a phenomenon applicable in almost all scientific fields when large numbers of single-issue measured phenomena are pooled together. Empirical evidences consistently show that most data sets of single-issue physical quantities such as time between earthquakes, earthquake's depth below the ground, population counts, quasar rotation rates, molar mass of frequently used chemical compounds, and so forth, are nearly perfectly logarithmic in their own right, individually considered, without being mixed with other data types, provided that order of magnitude of the bulk of the spread of the data is large enough (approximately 2.5 to 3.0 at a minimum).

When focus is shifted to the properties of density curves of logarithmic data sets instead of their digital behavior, it is found that such density curves must be generally falling in the aggregate, having a long tail to the right, and that the fall comes with a particular overall sharpness, all of which must also be well-coordinated between integral powers of ten. This changes the agenda, from looking at the way digits occur within typical numbers of everyday data, to looking at the typical forms of densities of typical pieces of everyday data; and it is proposed that the significant digits phenomenon (BL) is very much a phenomena of the propensity of everyday (positive) data to occur as lop-sided densities with tails to the right, falling off over the relevant intervals bounded by integral powers of ten with a variety of sharpness or rapidity, yet all having that same aggregate rate of fall over the entire range no matter. Such a particular fall coordinated and being in phase with those intervals bounded by integral powers of ten, results in that ubiquitous LOG(1+1/d) significant digits distribution.

Alternatively stated; it is not possible to separate leading digits distribution from the nature of the density curve itself. Leading digit distributions and density curves are not independent of each other, and any statement about leading digits configuration mathematically implies a restriction on the shape and range of the density curve. Consequently, this vista yields yet another essential perspective on BL regarding relative quantities, since a falling density curve on the right implies that small values are numerous and big values are few and far. Hence, as BL states a preference for lower digits, it also by implication states a preference for lower quantities as well.

Would it be possible to separate BL about digits from quantities all together? Can data be constructed in such a way that it obeys BL while possessing numerous big values and very few small values, namely an inverse quantitative configuration?

One possibility is to utilize a linear density rising from 9 to $2*10^8$ in a stop-and-go piecewise manner on the following 9 separate sections: {(9, 10), (80, 90), … , ($2*10^7$, $3*10^7$), ($1*10^8$, $2*10^8$)} deliberately constructed so that 1st order $LOG_{10}(1 + 1/d)$ is obeyed. But here all higher orders are inversely distributed and in violation of the general law. Else one could utilize mini k/x distributions over each of the same group of x-axis sections mentioned above where in general the big is more numerous than the small (overall) and BL is obeyed in the general sense. A third possibility is given by three discontinuous k/x curves superimposed, say, $k_1/x$ defined over (1, 10), $k_2/x$ defined over (10, 100), and $k_3/x$ defined over (100, 1000), where $k_3 > k_2$ and $k_2 > k_1$ such that the requirement $1/\ln(10) = k_1 + k_2 + k_3$ still holds, but curve on (100, 1000) is the highest, curve (1, 10) the lowest, and the big overtakes the small in the some general sense. Yet this quantitative reversal in lines with BL appears farfetched and extremely rare in real life typical data sets; and they are almost never empirically encountered. Almost all data sets show continuity and compactness; rarely broken down into separate and disjoint intervals. In the third scenario of 3 distinct k/x curves, the confluence of BL and continuity implies a consistent and steady quantitative decline, as all 3 curves must connect nicely at 10 and at 100 (in order to continue). All this is a strong confirmation that BL is very much a law about quantities, albeit indirectly so.

The consideration of (logarithmic) physical data sets, hand-made by mother nature herself who knows nothing about numbers and digits strongly suggests that the phenomenon is not relative to any number systems in use, but rather universal. Would it be possible to separate BL from our number system itself altogether or from any number system for that matter? The answer is in the affirmative. The law is scale invariant as well as base invariant. Moreover, and by far more significant, it is number-system invariant, reality floating out there in the universe, independent of whatever we do, calculate, count, invent, or define. To the extent that BL is a quantitative statement as well, this affords the law universality, rendering it a physical and scientific law.

Interestingly, Cartesian Coordinate System is number system invariant as well. The conceptual insight of the Cartesian Plane does not necessitate our number system, nor any other number system for that matter. It can be thought of as a purely geometrical construction, where distance and the primitive counting of unitary squares are playing the leading roles. Each unity of distance on the axes is duplicated and a count to the origin signifies quantity or position. Does it matter how we designate the units? No! A moment thought would convince everyone that data histograms and PDFs are number system invariant in essence (their pictorial aspect), and that the way the two axes are designated numerically is irrelevant quantitatively. Consequently, for a group of histograms representing a variety of physical (logarithmic) data sets, a singular and invariant quantitative BL measure should be true for all observers, regardless of number system in use or absence thereof [be they modern humans with base 10 positional number system, Babylonians, Romans, or Mayans.]

Clearly physics is number system invariant. The statement $GM_1M_2/R^2 = M_2A$ holds under base transformation. Moreover, and by far more significant, it is number system invariant, a primitive statement about quantities, not a numerical statement.

To further demonstrate the number system invariance property of BL, an imaginary snapshot of some logarithmic data sets relating to physical phenomenon is considered, a data set that can be visualized, say sizes of planets via the their 2-dimensional area being sensed. It's rather impossible to argue that this is all about digits or even about numbers, [and that absent a number system or digits no pattern or law could have been observed or stated]. On the contrary, it's all about relative quantities. Patterns in physical data sets transcend number systems and digits.

## [2] Bin Schemes

A primitive society without a number system whatsoever, a society that applies a single numerical symbol ('digit') for each and every numerical value (quantity), cannot observe nor formulate BL. Yet when such a society is presented with several histograms of logarithmic physical data sets, it should be able to find a common thread going through all of them, if not via LOG(1+1/d) since they have neither numbers nor digits, then via more direct (primitive) measure of relative occurrences of quantities. Yet, an attempt to simply cast a long system of repetitive bins along the x-axis, all with equal width throughout, and record relative values falling within each bin, cannot resolve this issue, no law can be observed; as all bins obtain almost equal portions of data, namely bin equality. This is so in spite of the fact that overall fall in all logarithmic histograms can be visually sensed! Surely a very crude bin width large enough to incorporate a significant portion of data does not show bin equality, but then results are unstable and depends on width. The fact that flat system cannot capture fall in histograms, is akin to higher orders of digits in BL which show near digital equality in spite of overall fall in densities.

What is meant by a bin system, or a bin scheme, is a repetitive positioning of equally-spaced imaginary D number of bins all along the x-axis, and where portion of data falling within each bin is recorded and finally aggregated. For example, a 3-bin system (D=3) can be imagined on (0, 2), (2, 4), (4, 6) as the first cycle, then on (6, 8), (8, 10), (10, 12) as the second cycle, and so forth. The width W here is 2 units per bin. The sub-intervals (0, 2) and (6, 8) belong to the same bin rank, namely the first bin. In general, flat (horizontal line) histograms, yield data with as many big quantities as small or medium ones. Rising histograms yield data with many more big quantities than small ones. Falling histograms yield many small ones and few big ones. Hence by constantly examining local relative fall (or rise) within bins and then aggregating all results, we construct a singular measure of relative quantities (i.e. overall fall/rise in histogram).

The term 'flat' signifies that the width of all the bins are not expanding but rather is made to be a constant in all cycles; that inflation factor F is 1 (as in the above 3-bin scheme having width W=2 throughout). Typically, the first cycle starts exactly from the origin 0, but it doesn't have to be always like that, result are extremely flexible.

Since so much of real life data is the result of multiplication processes, such as resultant data modeled on MCLT or exponential growth, values are 'stretched out' and 'expanded' along the x-axis 'rapidly' and 'forcefully' in a multiplicative manner, thus one needs to utilize an expanding bin scheme, letting the width of the bin expand multiplicatively by some inflation factor F > 1, in order to observe fall in histogram (relative quantities.) Clearly digital BL operates in the environment F = D + 1, as in all positional number systems. In this article, we allow for any formulation of bin systems, including F ≠ D + 1, as a pure measure of fall divorced from any connection to number systems.

Let us perform two expanding bin schemes. Scheme A has 4 bins with an inflation factor of 8 (meaning that after the first cycle of 4 bins, the width expands by a factor of 8, and continuing so for each new cycle, becoming ever lager). Scheme B has 7 bins with an inflation factor of 3. Naturally it is decided to perform both schemes A and B starting from 0. It is also deemed necessary to have the width starting quite narrowly at 0.0008 (for both schemes) in order to cast a refined net of bins applicable to data falling on (0, 1) as well. This aspect is similar to digital BL which operates under infinitely refined partitions on the left near the origin, and quite crude and wide partitions on the right. We shall use the notations: D for the number of bins ('digits'); F for the inflation factor; W for the initial bin width within the first set (cycle) of D bins on the left; S for the starting point of the whole scheme on the x-axis, namely the location of the left corner of the first bin in the first cycle (which is usually assigned as 0). The tables in Figures 1 and 2 show how logarithmic data sets and distributions all give almost the same bin spread for each of the two schemes. The summery of the two schemes' setups and results (average of all 9 data sets) are as follow:

Scheme A:
D=4 F=8 S=0 W=0.0008
Proportions = {**48.3%, 23.8, 15.9%, 12.0%**}

Scheme B:
D=7 F=3 S=0 W=0.0008
Proportions = {**22.6%, 18.4%, 15.2%, 13.1%, 11.3%, 10.1%, 9.3%**}

| Data Set | Bin A | Bin B | Bin C | Bin D |
|---|---|---|---|---|
| | | | | |
| Time Between Earthquakes | 48.3% | 25.0% | 15.3% | 11.5% |
| USA Population Centers | 48.9% | 23.1% | 16.0% | 12.0% |
| LOG Symmetrical Triangular (1, 3, 5) | 48.8% | 24.1% | 15.4% | 11.7% |
| k/x over (1, 1000000) | 49.3% | 21.7% | 16.3% | 12.7% |
| Exponential Growth, B=1.5, F=1.01 | 47.9% | 23.7% | 15.6% | 12.8% |
| Lognormal, Location=5, Shape=1 | 49.1% | 23.3% | 15.6% | 12.1% |
| Lognormal, Location=9.3, Shape=1.7 | 48.6% | 23.7% | 15.8% | 11.9% |
| Varied Data - Hill's Model | 46.3% | 25.3% | 16.2% | 12.2% |
| Chain U(U(U(U(U(0, 5666))))) | 47.8% | 24.1% | 16.1% | 12.0% |

**FIGURE 1** Four-bin Scheme A, D=4 F=8 S=0 W=0.0008

| Data Set | Bin A | Bin B | Bin C | Bin D | Bin E | Bin F | Bin G |
|---|---|---|---|---|---|---|---|
| | | | | | | | |
| Time Between Earthquakes | 22.4% | 18.1% | 15.5% | 13.1% | 11.7% | 10.1% | 9.0% |
| USA Population Centers | 22.6% | 18.9% | 15.4% | 13.1% | 10.9% | 9.9% | 9.1% |
| LOG Symmetrical Triangular (1, 3, 5) | 23.0% | 17.8% | 15.0% | 13.0% | 11.5% | 10.2% | 9.4% |
| k/x over (1, 1000000) | 21.5% | 17.9% | 15.2% | 13.4% | 12.5% | 10.3% | 9.2% |
| Exponential Growth, B=1.5, F=1.01 | 22.6% | 18.0% | 15.0% | 12.9% | 11.7% | 10.4% | 9.4% |
| Lognormal, Location=5, Shape=1 | 23.1% | 18.1% | 14.9% | 13.2% | 11.4% | 10.2% | 9.1% |
| Lognormal, Location=9.3, Shape=1.7 | 22.8% | 18.3% | 15.2% | 12.9% | 11.5% | 10.2% | 9.2% |
| Varied Data - Hill's Model | 22.0% | 20.1% | 15.6% | 13.1% | 10.3% | 9.5% | 9.4% |
| Chain U(U(U(U(U(0, 5666))))) | 23.2% | 17.8% | 15.6% | 13.3% | 10.8% | 10.2% | 9.1% |

**FIGURE 2** Seven-bin scheme B, D=7 F=3 S=0 W=0.0008

The data can be downloaded from The United States Geological Survey Organization website at http://earthquake.usgs.gov/earthquakes/eqarchives/epic/. The particular data set selected and obtained pertains to worldwide earthquakes, of any Richter scale magnitude, occurring at any depth below the ground, during the entire year of 2012. There were 19,452 earthquakes in 2012. The data set is 19,451 time intervals in seconds between earthquake occurrences.

USA Census data on population counts of all incorporated cities and towns relates to 19,509 such population centers in the year 2009. This data set adheres to BL very closely. The data can be downloaded from the website: http://www.census.gov/popest/data/historical/2000s/vintage_2009/datasets.html with the choice of “Vintage 2009 City and Town (Incorporated Place and Minor Civil Division) Population Datasets”; the selection of “All States” as the file; and the choice of the last column titled “POP_2009” as the data set.

The Log symmetrical triangular data set relates to the density of LOG of simulated data, which starts from log = 1, centers on log =3, ends on log = 5, and it gives rise to data by simply calculating $10^{\text{Triangular}}$. This data set is perfectly logarithmic.

The exponential growth data set represents the first 10,000 elements from exponential 1% growth from a base of 1.5.

'Varied Data' refers to 34,269 values randomly obtained from 70 different Internet sources and topics. Such random selection is (approx) logarithmic as per the mathematically rigorous proof in BL about distribution of infinitely many distributions.

Chain of distributions refers to the Uniform(a, b) distribution, where parameter a is always 0, while b is not fixed/constant as is the standard, but rather itself a random variable derived from yet another Uniform distribution, etc.
In symbols: U(0, U(0, U(0, U(0, U(0, 5666))))). The last/inner Uniform(0, 5666) represents a standard distribution without parametrical dependency on any other distribution. This chain is nearly perfectly logarithmic, and precisely so in the limit as the length (number of sequences) of such parametrical dependency goes to infinity.

There are similarities and differences between partitioning of the x-axis by way of bin systems and by way of positional number systems along digital lines as in BL. Positional number systems contain infinitely many digital cycles which are not 'directly' or 'intentionally' placed from the 0 origin onwards as can be done actively in bin systems. Rather positional number systems can be thought of as creating digital cycles from the base onwards, as well as from the base backwards towards the origin in an infinite process, which in the limit 'starts' from the origin.

The tables in Figures 1 and 2 clearly tell the story of an overall definite fall in the densities of all of these data sets and distributions, no matter how many bins are set to measure it or the exact value of F. In addition, for a given number of D bins, F expansion factor, S and W [i.e. for a given bin scheme], the rate of fall as measure by the relative proportions of overall data falling within each bin is remarkably consistent across all (logarithmic) data types and distributions, with minor variations.

## [3] Non-Expanding Bin System Measuring Fall in k/x Distribution

Postulating that the generic pattern in how relative quantities are found in nature is such that the frequency of quantitative occurrences is inversely proportional to quantity, we are then led to explore results from bin systems fitting k/x distribution, the density curve whose arrangement of relative quantities constitutes exactly such a pattern.

Three cases shall be examined: non-expanding, once-expanding, and twice-expanding bin systems. The goal is to arrive at a general expression for bin proportions in the k/x case, which would then be checked against empirical bin results of real data - assuming that a consistent law could be mathematically found for the k/x case.

We shall evaluate a single cycle D-bin scheme on the density curve k/x distribution defined over (w, (D+1)w) with bins standing between (w, 2w), [2w, 3w), … , [Dw, (D+1)w). This represents a basic or unitary bin system with no expansion, having only one cycle. Five features are involved in this construction:

(I) Avoidance of an upward explosion start of the density at the origin 0 which would be undefined due to a division by 0.
(II) Equal spacing (width) of all D bins.
(III) Equality between the width of the bins and the separation of launch from the 0 origin, namely that the w length of the step from the origin to the launch of k/x is also the length of the first bin and all subsequent bins, symbolically written as (2w - w) = (w - 0).
(IV) No coordination is employed or attempted whatsoever with any number system or significant digits on the x-axis below.
(V) No assumption is made about the value of the width w, which is left completely flexible to take on any value whatsoever, even fractional value.

Equating the entire area to one, we obtained ∫ k/x dx = 1 over [w, (D+1)w)], therefore k[ln((D+1)w) - ln(w)] = 1, or k[ln(D+1) + ln(w) - ln(w)] = 1, so that k[ln(D+1)] = 1, hence k = 1/ln(D+1). Evaluating the portion of area hanging over bin #d (with d running from 1 to D, as in digits), we obtain P(d) = ∫ k/x dx over [d*w, (d+1)*w)], namely P(d) = [1/ln(D+1)]*[ln(d+1) + ln(w) – ln(d) – ln(w)] = [1/ln(D+1)]*[ln(d+1) – ln(d)], or P(d) = [1/ln(D+1)]*ln[(d+1)/(d)], and finally **P(d) = ln(1+1/d) / ln(D+1)**.

Double applications the logarithmic identity $LOG_A X = LOG_B X / LOG_B A$ yield [log(1+1/d) / log(e)] / [log(D+1) / log(e) ], and finally **P(d) = log(1+1/d) / log(D+1)**. Hence the above bin probability expression is perfectly compatible with the more general BL encompassing other bases, since D + 1 is the equivalent value of the base in any number system. Yet, this whole scheme must be considered with severe reservation due to its intrinsic limitation in lacking expansions of the cycles. Without expansions to ever larger bin cycles, consistency of results is in doubt.

Nonetheless, the remarkable philosophical and conceptual implication of this result is that the form of BL in the expression log(1+1/d)/log(base) serves as a general quantitative proportional law outside any digital framework (in some restricted cases, namely k/x without bin expansion, plus others). This is so since width w cancels out and drops from the calculations, as a result w could take on any value without affecting the above expression. The significance of all this emanates from the fact that a variety of 1st digits might be mixed within any single bin! Also, the fact that results are independent of the width w is significant, since this lends the bin scheme universality and consistency. Equivalently stated, log(1+1/d)/log(base) serves also to express relative quantities and the fall in the density of k/x defined over some restricted range. As an example, data set corresponding perfectly to k/x defined over (20, 200) is analyzed as a non-expanding bin scheme with w = 20 and D = 9. The 9 bins measuring the fall are {(20, 40), (40, 60), (60, 80), (80, 100), (100, 120), (120, 140), (140, 160), (160, 180), (180, 200)}. Surely 30.1% of data falls within the first bin (20, 40); 17.6% of data falls within (40, 60); while

only 4.6% falls within the last bin (180, 200); and all proportions perfectly fit the expression $LOG_{10}(1 + 1/d)$ as in BL. Yet the two first significant digits 2 and 3 which normally compete for leadership now reside harmoniously within bin (20, 40); the two adversarial first digits 4 and 5 now live peacefully within bin (40, 60); and so forth.

It should be noted that the above setup of k/x distribution defined over [w, (D+1)w)] is considered perfectly logarithmic in digital BL sense for any base B, where D - the numbers of bins - signifies the number of all first digit possibilities (such as 9 first digits in base 10). This is so since exponent difference of this range is $LOG_B((D+1)w) - LOG_B(w) = LOG_B(D+1)$, and since D+1 represents the base as well, this reduces to $LOG_B(B)$, namely 1, which is an integral number.

## [4] Once-Expanding Bin System Measuring Fall in k/x Distribution

Expanding bin scheme only once for the density k/x over (w, (D+1)w + DFw) means that the original w width of the bin is being inflated by F on the second cycle. In other words, where the width of each bin in the 2nd cycle is not w but Fw.

1st cycle is on (w, 2w), [2w, 3w), … , [Dw, (D+1)w).
2nd cycle is on [(D+1)w, (D+1)w + Fw), [(D+1)w + Fw, (D+1)w + 2Fw), … ,
, … , [(D+1)w + (D-1)Fw, (D+1)w + DFw).

Equating the entire area to one, we obtained ∫ k/x dx = 1 over [w, (D+1)w + DFw], hence k[ ln(w*[(D+1) + (DF)]) - ln(w)] = 1, or k[ ln(w) + ln[(D+1) + (DF)] - ln(w)] = 1, so that k[ ln[(D+1) + (DF)] ] = 1, and finally k = 1/ ln(1 + D + DF). Evaluating the **first** portion of area (1st cycle) hanging over bin #d (d running from 1 to D as in digits), we obtain:

P1(d) = ∫ k/x dx over [d*w, (d+1)*w)]
P1(d) = [1/ln(1 + D + DF)]*[ln(d+1) + ln(w) – ln(d) – ln(w)]
P1(d) = [1/ln(1 + D + DF)]*[ln(d+1) – ln(d)].

Evaluating the **second** portion of area (2nd cycle) hanging over bin #d (with d running from 1 to D, as in digits), we obtain:

P2(d) = ∫ k/x dx over [(D+1)w + (d-1)*Fw, (D+1)w +(d)*Fw],
P2(d) = [1/ln(1 + D + DF)]*[ln((D+1)+dF) + ln(w) – ln((D+1)+(d-1)F) – ln(w)],
P2(d) = [1/ln(1 + D + DF)]*[ln((D+1)+dF) – ln((D+1)+(d-1)F)]

Now combining both areas, namely P(d) = P1(d) + P2(d), we finally get:

P(d) = [1/ln(1 + D + DF)]*[ ln(d+1) – ln(d) + ln((D+1)+dF) – ln((D+1)+(d-1)F)]

And simplifying by using the identity LOG(A) - LOG(B) = LOG(A/B) we finally get:

P(d) = [ ln(1 +1/d) + ln[ (1+D+dF)/(1+D+(d-1)F) ] / [ ln(1 + D + DF) ]

## [5] Once-Expanding Bins for k/x Reduces to BL whenever F = D + 1

The Bin system set up for the k/x distribution with one expansion cycle reduces to BL whenever inflation factor F is made equal to (D+1) as in all proper number systems. In other words, whenever F = D + 1, that once-expanding bin system for k/x reduces to non-expanding bin system for k/x. The implication here is quite significant, because it means that the very act of (a singular) expansion (doubling the number of bins) does not change bin proportions in any way in the k/x case whenever F = D + 1.

To prove the assertion, (D + 1) is simply substituted for F into the above P(d) expression, hence:

P(d) = [ ln(1 +1/d) + ln[ (1+D+dF)/(1+D+(d-1)F) ] / [ ln(1 + D + DF) ]
P(d) = [ ln(1 +1/d) + ln[ (1+D+d(D + 1))/(1+D+(d-1)(D + 1)) ] / [ ln(1 + D + D(D + 1))]
P(d) = [ ln(1 +1/d) + ln[ ((1+D)(1 + d))/((1+D)(1+(d-1)) ] / [ ln(1 + D + $D^2$ + D))]
P(d) = [ ln(1 +1/d) + ln[ ((1+D)(1 + d))/((1+D)(d)) ] / [ ln(1 + 2D + $D^2$))]
P(d) = [ ln(1 +1/d) + ln[ (1 + d)/(d) ] / [ ln(1 + 2D + $D^2$))]
P(d) = [ ln(1 +1/d) + ln[ (1 + 1/d) ] / [ ln($(D + 1)^2$)]
P(d) = 2*ln(1 +1/d) / ln($(D + 1)^2$)
P(d) = 2*ln(1 +1/d) / 2*ln(D + 1)
P(d) = ln(1 +1/d) / ln(D + 1)
P(d) = log(1 +1/d) / log(D + 1)
P(d) = log(1 +1/d) / log(Base)
P(d) = BL

## [6] Twice-Expanding Bin System Measuring Fall in k/x Distribution

Expanding bin scheme twice for the density k/x over (w, (D+1)w + DFw + $DF^2$w) means that the original bin's width w is being inflated by F on the second cycle, and then inflated by F-squared on the third cycle.

1st cycle is on (w, 2w), [2w, 3w), … , [Dw, (D+1)w).
2nd cycle is on [(D+1)w, (D+1)w + Fw), [(D+1)w + Fw, (D+1)w + 2Fw), … ,
, … , [(D+1)w + (D-1)Fw, (D+1)w + DFw).
3rd cycle is on [(D+1)w + DFw, (D+1)w + DFw + FFw), … , … ,
, … , [(D+1)w + DFw + (D-1)FFw), (D+1)w + DFw + DFFw).

Equating the entire area to one, we obtained ∫ k/x dx = 1 over the range of [w, (D+1)w + DFw + $DF^2$w)], hence k[ ln(w) + ln[(D+1) + DF + $DF^2$] - ln(w)] = 1, or k[ ln[(D+1) + DF + $DF^2$] ] = 1, and finally k = 1/ ln(1 + D + DF + $DF^2$).

Evaluating the **first and second** portion of areas yields the same results as in the once-expanding bin system except for the k constant which differs here.
Evaluating the **third** portion of area hanging over bin #d (with d running from 1 to D, as in digits), we obtain:

$P3(d) = \int k/x \, dx$ over $[(D+1)w + DFw + (d-1)*F^2w, \ (D+1)w + DFw + (d)*F^2w]$
$P3(d) = k ( \ln[w] + \ln[(D+1) + DF + (d)*F^2] - \ln[w] - \ln[(D+1) + DF + (d-1)*F^2] )$
$P3(d) = k ( \ln[(D+1) + DF + (d)*F^2] - \ln[(D+1) + DF + (d-1)*F^2] )$
$P3(d) = k* \ln( [(D+1) + DF + (d)*F^2] / [(D+1) + DF + (d-1)*F^2] )$

Combining all 3 areas, namely P(d) = P1(d) + P2(d) + P3(d), we finally get:

$P(d) = k*\mathbf{ln}(1 +1/d) + k*\mathbf{ln}([1+D+dF]/[1+D+(d-1)F]) +$
$\quad k*\mathbf{ln}([(D+1) + DF + (d)*F^2]/[(D+1) + DF + (d-1)*F^2])$

$$\mathbf{P(d)} = \frac{\ln(1 + 1/d) + \ln\left(\frac{[\,1 + D + (d)F\,]}{[\,1 + D + (d-1)F\,]}\right) + \ln\left(\frac{[\,1 + D + DF + (d)*F^2\,]}{[\,1 + D + DF + (d-1)*F^2\,]}\right)}{\ln(1 + D + DF + DF^2)}$$

## [7] Twice-Expanding Bins for k/x Reduces to BL when F = D + 1

The twice-expanding bin system set up for k/x distribution also reduces to BL whenever inflation factor F is made equal to (D+1) as in all proper number systems. In other words, whenever F = D + 1, that twice-expanding bin system reduces to non-expanding bin system. The implication here is that the very act of doubly expanding the bins does not change bin proportions in any way whenever F = D +1.

To prove the assertion, (D + 1) is simply substituted for F into the above P(d) expression, hence:

$$\frac{\ln(1+1/d) + \ln\left(\frac{[\,1 + D + (d)(D+1)\,]}{[\,1 + D + (d-1)(D+1)\,]}\right) + \ln\left(\frac{[\,1 + D + D(D+1) + (d)*(D+1)^2\,]}{[\,1 + D + D(D+1) + (d-1)*(D+1)^2]}\right)}{\ln(1 + D + D(D+1) + D(D+1)^2)}$$

$$\frac{\ln(1+1/d) + \ln\left(\frac{[\,(d+1)(D+1)\,]}{[\,(d - 1 + 1)(D+1)\,]}\right) + \ln\left(\frac{[\,(1+D)(D+1) + (d)*(D+1)^2\,]}{[\,(1+D)(D+1) + (d-1)*(D+1)^2]}\right)}{\ln(\,(1+ D)(D+1) + D(D+1)^2)}$$

$$\frac{\ln(1+1/d) + \ln\left(\frac{[\,(d+1)(D+1)\,]}{[\,(d)(D+1)\,]}\right) + \ln\left(\frac{[\,(1+D)*[\,(D+1) + (d)*(D+1)\,]]}{[\,(1+D)*[(D+1) + (d-1)*(D+1)\,]]}\right)}{\ln(\,(1+D)*[(D+1) + D(D+1)]\,)}$$

$$\frac{\ln(1+1/d) + \ln\left(\frac{[d+1]}{[d]}\right) + \ln\left(\frac{[\,(D+1) + (d)*(D+1)\,]}{[(D+1) + (d-1)*(D+1)\,]}\right)}{\ln(\,(1+D)*[(D+1)(1+D)]\,)}$$

$$\frac{\ln\left(\frac{[d+1]}{[d]}\right) + \ln\left(\frac{[d+1]}{[d]}\right) + \ln\left(\frac{[(D+1)[1+d]]}{[(D+1)[1+(d-1)]\,]}\right)}{\ln((1+D)^3)} = \frac{3*\ln\left(\frac{[1+d]}{[d]}\right)}{3*\ln(1+D)} = \mathbf{BL}$$

Another check on a triple-expanding bin system for k/x distribution was perform for F = D + 1, confirming the same reduction to BL as seen in the previous two cases. While this does not constitute a formal proof, it strongly suggests that any N-expanding (and also ∞-expanding) bin systems constructed for k/x would yield the same reduction to the BL given that F = D + 1.

## [8] Infinitely-Expanding Bin System Measuring Fall in k/x A.K.A. The General Law of Relative Quantities

Algebraic expressions for bin proportions of k/x distribution for higher expansion orders perfectly follow the above pattern as a sequence of ever increasing terms in the numerator and in the denominator. The first 4 elements of this sequence, beginning with a non-expanding bin system, and ending with a triple-expanding bin system are as follow:

$$\frac{\ln\left(\dfrac{[\,1 + (d)\,]}{[\,1 + (d-1)\,]}\right)}{\ln(1 + D)}$$

$$\frac{\ln\left(\dfrac{[\,1 + (d)\,]}{[\,1 + (d-1)\,]}\right) + \ln\left(\dfrac{[\,1 + D + (d)F\,]}{[\,1 + D + (d-1)F\,]}\right)}{\ln(1 + D + DF)}$$

$$\frac{\ln\left(\dfrac{[\,1 + (d)\,]}{[\,1 + (d-1)\,]}\right) + \ln\left(\dfrac{[\,1 + D + (d)F\,]}{[\,1 + D + (d-1)F\,]}\right) + \ln\left(\dfrac{[\,1 + D + DF + (d)F^2\,]}{[\,1 + D + DF + (d-1)F^2\,]}\right)}{\ln(1 + D + DF + DF^2)}$$

$$\frac{\ln\left(\dfrac{[1+(d)]}{[1+(d-1)]}\right) + \ln\left(\dfrac{[1+D+(d)F]}{[1+D+(d-1)F]}\right) + \ln\left(\dfrac{[1+D+DF+(d)F^2]}{[1+D+DF+(d-1)F^2]}\right) + \ln\left(\dfrac{[1+D+DF+DF^2 + (d)F^3]}{[1+D+DF+DF^2+(d-1)F^3]}\right)}{\ln(1 + D + DF + DF^2 + DF^3)}$$

Convergence of this long sequence of bin proportions of the k/x distribution over an infinite expansion has been suggested by way of finite computer simulations (calculations) for a variety of F and D values, although of course this does not prove convergence.

One important aspect in all expanding bin systems is independence of k/x results on initial width w, as it drops out in the calculations and thus could take on any value without affecting resultant expression, just as was seen in the non-expanding case. Yet caution should be exercised when real life data is concerned and where large values of w totally distort results. The dichotomy emanates from the fact that for all bin schemes on k/x, by definition there exist no values/data between 0 and w, while for real life data sets huge portion of data might be hiding there. In other words, k/x schemes start at w, while real data may start at or near 0. Therefore, setting w equal to some very small fractional value is advisable if compatibility and correspondence between the theoretical framework developed here about the generic k/x distribution and real life data sets is desired. Even having the smallest value (min) in any given data set larger than w is not a guarantee that w was chosen small enough if lots of values in the data set congregate just above min, hence w should be made extremely small in relation to all possible data sets under consideration.

Note: All the mathematical derivations on the relative quantities of k/x distribution assume equality between the width of the bins and the separation from the origin 0, namely that w - the length of the step from the origin to the launch of k/x - is also the width of the bins in the first cycle. Violation of this equality leads to different results. If separation is longer than initial bin sizes, then greater bin equality prevails. If separation is shorter than bin sizes, then sharper fall is observed and extreme bin inequality prevails (severe bin skew-ness that is).

## [9] Confirmation Matching k/x Fall with Empirical Bins on Real Data

Let us check real life data sets and abstract logarithmic distributions all viewed through bin scheme's prism, against the (apparent) converging limit of the algebraic sequence derived here regarding the fall between the bins of the generic density k/x infinitely expanded (utilizing corresponding bin systems, namely bin schemes having the same values of F and D in both theoretical and empirical cases for valid comparisons and confirmation).

The two bin schemes in Figures 1 and 2 are compared with equivalent infinitely-expanded bins superimposed on k/x [using computer calculation of ever increasing cycles of expansions on k/x, leading to apparent convergence]. It shows a near perfect agreement between empirical and theoretical results. The two averages of the bin results from all nine data sets and distributions, together with equivalent k/x bin results are as follow:

Scheme A:
D=4  F=8  S=0  W=0.0008
Average proportions of all 9 real data sets & dist  = {48.3%, 23.8%, 15.9%, 12.0%}
Limit of algebraic seq. of k/x dist ∞-expanded       = {48.6%, 23.7%, 15.8%, 11.9%}

Scheme B:
D=7  F=3  S=0  W=0.0008
Avg of all 9 real data sets & dist = {22.6%, 18.4%, 15.2%, 13.1%, 11.3%, 10.1%, 9.3%}
Limit of seq. k/x ∞-expanded       = {22.9%, 18.3%, 15.2%, 13.0%, 11.4%, 10.1%, 9.1%}

The limit of the infinite sequence of algebraic expressions for the infinitely expanded k/x distribution should now be construed as constituting the general law governing occurrences of relative quantities in all logarithmic data sets (as a function of the manner in which such a concept is measured, namely as a function of D and F, with S and w being insignificant and irrelevant as they are confined to very small fractional values and do not affect results.)

## [10] 9-Bin Systems with F=10 on Real Data All Yield $LOG_{TEN}$(1+1/d)

As an additional confirmation of the general bin theory developed here, the same group of 9 logarithmic data sets and distributions shall be examined under a 9-bin system having inflation factor F = 10. Results are then compared to the logarithmic distribution $LOG_{10}$(1+1/d) and are found to fit nicely. Although such a bin system may appear identical to our numbers system, digits, and BL, as if imitating them, this is clearly not the case for two reasons, (I) such a bin system may start at any point including 0 (as long as it's not too far from the origin) and still yield the same results, (II) the width of the first bin has a finite value and it may be of a substantial size (within a limit, not too large) still yielding the same results. Consequently, and most significantly, the bins are not at all aligned and coordinated on the digital postmarks of our number system (each bin within each cycle contains a variety of significant 1st digits mixed in!) What this results show is that these digital postmarks are not the only way to measure relative quantities, and their absence here has no effect on results at all, which are $LOG_{10}$(1+1/d) just the same. This should be considered as a decisive argument that $LOG_{10}$(1+1/d) is for the most part all about relative quantities, and its digital application is but a minor event in the much larger quantitative drama.

The table in Figure 3 depicts results from the same group of 9 real life logarithmic data sets and distributions viewed through the lens of a 9-bin system with F=10 starting at 0.033 and having an initial width 0.07. The fit into $LOG_{10}(1+1/d)$ is quite satisfactory! That very refine and narrow width start of 0.07, and the positioning of the beginning of the whole bin scheme quite near the origin at 0.033 were two essential features contributing to the 'success' of the results in term of closeness to the logarithmic. The table in Figure 4 depicts quite different results from the same data sets and distributions viewed through a much cruder lens of 9-bin system with F=10 starting too late at 5, and having an initial very thick width of 311. This shows that bin systems need starting near origin with some very small and refined initial width.

Here we truly encounter BL and the logarithmic (30.1%, 17.6%, … , 4.6%) in its most general form, without digits, free and independent of any number system whatsoever, and this is so for practically all real life data sets and abstract distributions that are known to be 'logarithmic'. Crucially, such LOG(1 + 1/d) bin results are theoretically supported by the corresponding bin schemes on the generic k/x distribution as demonstrated earlier (where F = D + 1, conjectured to reduce to BL no matter how many cycles are used.)

| Data Set | Bin A | Bin B | Bin C | Bin D | Bin E | Bin F | Bin G | Bin H | Bin I |
|---|---|---|---|---|---|---|---|---|---|
| | | | | | | | | | |
| Time Between Earthquakes | 29.2% | 17.8% | 13.3% | 10.5% | 8.1% | 6.6% | 5.6% | 4.8% | 4.2% |
| USA Population Centers | 30.3% | 17.3% | 12.8% | 9.7% | 7.3% | 6.8% | 5.7% | 5.3% | 4.7% |
| LOG Sym. Triangular (1, 3, 5) | 29.6% | 17.9% | 12.5% | 9.6% | 8.0% | 6.8% | 5.8% | 5.0% | 4.8% |
| k/x over (1, 1000000) | 30.2% | 17.6% | 12.3% | 9.9% | 7.8% | 6.8% | 5.9% | 5.0% | 4.5% |
| Exp. Growth, B=1.5, F=1.01 | 27.4% | 18.0% | 13.0% | 10.1% | 8.2% | 7.0% | 6.0% | 5.4% | 4.8% |
| Lognormal, Loc=5, Shape=1 | 30.7% | 18.4% | 12.6% | 9.3% | 7.6% | 6.3% | 5.6% | 5.0% | 4.5% |
| Lognormal, Loc=9.3, Shape=1.7 | 30.4% | 17.4% | 12.4% | 9.6% | 7.9% | 6.8% | 5.8% | 5.2% | 4.5% |
| Varied Data - Hill's Model | 32.4% | 16.4% | 11.2% | 8.9% | 8.6% | 6.1% | 6.3% | 5.5% | 4.6% |
| Chain U(U(U(U(U(0, 5666))))) | 30.7% | 17.3% | 12.0% | 9.6% | 8.1% | 6.4% | 6.0% | 5.6% | 4.4% |
| General R.Q. Law D=9 F=10 | 30.1% | 17.6% | 12.5% | 9.7% | 7.9% | 6.7% | 5.8% | 5.1% | 4.6% |

**FIGURE 3** 9-Bin Scheme F=10 Yields ≈ $LOG_{10}(1+1/d)$ (Start=0.033 Width=0.07)

| Data Set | Bin A | Bin B | Bin C | Bin D | Bin E | Bin F | Bin G | Bin H | Bin I |
|---|---|---|---|---|---|---|---|---|---|
| | | | | | | | | | |
| Time Between Earthquakes | 35.9% | 17.7% | 12.3% | 9.0% | 7.5% | 5.6% | 4.7% | 3.8% | 3.4% |
| USA Population Centers | 36.1% | 21.3% | 12.9% | 8.8% | 6.1% | 4.8% | 4.1% | 3.3% | 2.6% |
| LOG Sym. Triangular (1, 3, 5) | 42.0% | 18.4% | 11.3% | 8.0% | 6.1% | 4.7% | 3.8% | 3.1% | 2.6% |
| k/x over (1, 1000000) | 51.3% | 15.8% | 8.8% | 6.1% | 4.8% | 4.3% | 2.8% | 3.2% | 2.9% |
| Exp. Growth, B=1.5, F=1.01 | 39.3% | 15.8% | 11.0% | 8.4% | 6.8% | 5.7% | 4.9% | 4.4% | 3.8% |
| Lognormal, Loc=5, Shape=1 | 77.6% | 14.9% | 4.3% | 1.7% | 0.8% | 0.4% | 0.2% | 0.1% | 0.1% |
| Lognormal, Loc=9.3, Shape=1.7 | 29.9% | 18.2% | 12.7% | 9.8% | 7.6% | 6.7% | 5.8% | 4.9% | 4.4% |
| Varied Data - Hill's Model | 53.2% | 13.8% | 9.4% | 6.1% | 4.7% | 3.6% | 4.0% | 3.0% | 2.3% |
| Chain U(U(U(U(U(0, 5666))))) | 80.3% | 11.4% | 4.3% | 1.8% | 1.0% | 0.6% | 0.3% | 0.2% | 0.1% |
| General R.Q. Law D=9 F=10 | 30.1% | 17.6% | 12.5% | 9.7% | 7.9% | 6.7% | 5.8% | 5.1% | 4.6% |

**FIGURE 4** 9-Bins F=10 Greatly Deviates from $LOG_{10}(1+1/d)$ (Start=5 Width=311)

One crucial aspect of all bin schemes and bin laws is that empirically exact values of starting point S and initial width w do **not** matter at all (almost), as long as S and w are made quite small. In other words, bin results are not really affected when changes are made to S and w. Only when values of S and w become quite large that bin proportions became increasingly more and more dependent on them. For low values of S and w, the only factors in bin results are: (I) number of bins D and (II) inflation factor F. In fact, moderate changes in S do not affect bin proportions in the least, and that is true for logarithmic as well as for non-logarithmic data. Moderate changes in w on the other hand yield very mild fluctuations in results for logarithmic data, and huge/dramatic fluctuations for non-logarithmic data, It is advisable in any case to universally stick to S = 0 and w < 0.001 in **all** bin schemes. An arbitrarily imposed **standardization rule** such as S = 0 and w = 0.0005 for **all** bin schemes could facilitate consistent comparisons in the field [unless data set under consideration contains some values on (0, 0.0005), in which case w should be made even smaller than 0.0005.] In order to formally define the bin model mathematically, for data sets with only non-zero positive numbers, the starting point S is made to be at the origin, and the initial bin width w is made infinitesimally small by defining it as a limiting process where w approaches zero.

## [11] Correspondence in Data Classification between Bin Systems & BL

An extremely crucial fact that must be acknowledged here is that [for a data set] being digit-wise logarithmic in the context of BL, and obeying in general all bin systems laws, go together. Having non-logarithmic digit configuration corresponds to having a bin configuration which differs than the one gotten for logarithmic data sets. Yet this fact should not be surprising in the least, because digital BL is nothing but one particular example of the generic idea of the bin scheme, a singular manifestation of a much larger universe of bin system possibilities. A special one to be sure, due to F = D + 1, and due to width w being made infinitesimally small towards the origin, but a bin scheme nonetheless.

In conclusion, being 'logarithmic' or 'non-logarithmic' is an absolute and universal property of any data set, irrespective of the base or bin scheme in use [in other words: it's a quantitative property, not numerical, not digital]. Curiosity compels one to ask then: "What is the basic or the most essential characteristic of being 'logarithmic'?" An appealing answer to this dilemma is the statement that logarithmic-ness is the intrinsic property of having an overall (average) decrease in relative quantities in the same rapidity and manner as that of the k/x generic distribution. Perhaps one may sum up uniqueness and simplicity in the fall of k/x by pointing to the fact that by <u>doubling</u> x value we cut the density exactly by a <u>half</u> (hence sum invariance characterization principle.) For example, for 0.4342945/x over (10, 100), density height or histogram count on x = 40 is exactly half that on x = 20. Such an exact relationship between x and its quantitative frequency is of course unique to k/x distribution.

It must be emphasized that such a definition of 'logarithmic-ness' utilizes k/x defined over a 'truly long range', namely k/x over (≈ a, ≈ ∞). For all practical purposes one can generate k/x data on say (1, $10^{10}$) and call it logarithmic. But k/x on (1, 10) is not logarithmic in this sense. It only complies with digital BL when base 10 is used, not when other bases are used, and it certainty does not comply with any bin laws unless it exactly mimics BL base 10. Short-ranged k/x are 'too sensitive' to base or bin changes. Long-ranged k/x are crude and stable, being base and bin scheme invariant, not in the sense that proportions are the same in different bases or bins, there are not, but in the sense that long-range k/x obeys all digital and bin laws regardless. The main challenge for k/x defined over a short range in obeying digit and bin laws is a perfect match with the measuring cycle, namely that k/x does not abruptly start or terminate in the middle of a bin or a digit cycle which strongly discriminates against some bins/digits, but rather that it starts and terminates at the beginning/end of a cycle. For k/x defined over a very long range, this issue is a minor one, since any such abrupt end or launch incorporating a partial cycle has a very small effect in the grand scheme of things. All this is also nicely consistent with the series of algebraic expressions for k/x where $F \neq D + 1$, since convergence to the limit necessitates the consideration of a 'very long range' on the x-axis in an infinitely-expanded bin system. All this is also nicely consistent with the consideration of the General Law representing the limit of the sequence of algebraic expressions of bin systems for k/x where $F \neq D + 1$, since convergence to the limit necessitates the consideration of a 'very long range' on the x-axis in an infinitely-expanding cycles of bins! Approximate convergence there is not found in merely 5 or 10 cycles say, but rather it is often obtained after say 30 or 100 full cycles, while for some D and F combinations it may take hundreds or longer just to get somewhat close enough to the true limit. Any real positive data terminates at some 'very finite' P point, while the grid of the bin system measuring (fall in) relative quantities within the given data constitutes an **infinite** net of bins emerging from the **left** towards P. The General Law on the other hand is derived from **infinitely** expanding bins **rightwards** from w to ∞ for the k/x distribution. Yet in spite of this '**directional dichotomy'** both, bin proportions of real finite data, and the General Law obtained from k/x, perfectly correspond. As noted in chapter 3, comment III, all the mathematical derivations for the General Law on the relative quantities of k/x distribution assume equality between the width of the bins in the first cycle and the separation from the 0 origin at launching. Violation of this equality leads to different results in the finite expansion case, where separation that is longer than the size of the first bin yields greater bin equality, and separation that is shorter than the size of the first bin yields sharper fall and extreme bin inequality. However, in the limit when the number of expansions goes to infinity this crucial constraint becomes irrelevant, and results are independent of the exact launching point. This fact provides us with one more crucial brick in the whole harmonious and consistent bin edifice, guaranteeing that any bin result of any particular k/x defined over a very large interval (order-of-magnitude-wise), would closely correspond to the theoretical limit of the bin proportions of infinite cycles construction built onto the generic k/x distribution [The General law].

It must be noted though that **all** random data sets do not possess (directly) that k/x property of having direct proportion to 1/x in density, nor is there an exact halving in density whenever x is doubled (hence all random data sets do **not** agree with that sum invariance characterization principle in BL). Yet random data relate to k/x indirectly by way of having corresponding overall fall in density (in the aggregate) as measured via bin schemes. All random data shows graduation and development in fall as shall be discussed in section 14, and their LOG density appears Normal-like or as an upside-down-U-like shaped curve. Only LOG density of the k/x distribution is uniform and steady throughout.

Base invariance principle in BL can now be interpreted as the principle of the universality of the 'logarithmic-ness' property. That no matter what bin scheme is used, no matter what base is applied in digital BL, classification of a given data set is a constant and universal, namely measuring-system-invariant. [That a change in B base, F, or D, does not revolutionize data classification.]

Let us give a more concrete statement of the universality of data classification as follows: While any given logarithmic data set consistently and simultaneously obeys all possible bin laws [the General Law] given <u>any</u> F or D value; any given non-logarithmic data set is a serial violator of the law of relative quantities, simultaneously disobeying all possible D and F versions of it.

## [12] The Remarkable Malleability of Bin Schemes

Surprisingly, we may vary the value of inflation factor F within a single bin scheme and still get a consistent pattern (law) across all logarithmic data sets. A 5-bin scheme, starting at 0, with an initial width 0.007, having an arbitrary and finite inflation factor vector $F_i$ = {2, 3, 4, 2, 5, 3, 6, 3, 5, 7, 4, 2, 3, 2, 7, 8, 9, 7, 3, 6} yields consistent results as seen in Figure 5. Although expansion along the x-axis is normally achieved by way of infinitely applying a fixed inflation factor, here just the width of the last bin cycle is sufficiency large to enclose the entire range of each data set and distribution, since (5)*(0.007)*2*3*4*2*5*3*6*3*5*7*4*2*3*2*7*8*9*7*3*6 > Max of Each Data Set.

| Data Set | Bin A | Bin B | Bin C | Bin D | Bin E |
|---|---|---|---|---|---|
| Time Between Earthquakes | **38.7%** | **21.8%** | **15.8%** | **12.6%** | **11.1%** |
| USA Population Centers | **36.0%** | **22.7%** | **16.5%** | **13.6%** | **11.2%** |
| LOG Symmetrical Triangular (1, 3, 5) | **36.6%** | **22.7%** | **16.7%** | **13.3%** | **10.7%** |
| k/x over (1, 1000000) | **34.2%** | **22.7%** | **17.4%** | **13.3%** | **12.4%** |
| Exponential Growth, B=1.5, F=1.01 | **35.3%** | **23.0%** | **16.9%** | **13.5%** | **11.3%** |
| Lognormal, Location=5, Shape=1 | **33.6%** | **23.3%** | **17.7%** | **13.9%** | **11.4%** |
| Lognormal, Location=9.3, Shape=1.7 | **35.5%** | **23.2%** | **17.0%** | **13.3%** | **10.9%** |
| Varied Data - Hill's Model | **35.6%** | **22.2%** | **16.8%** | **12.8%** | **12.6%** |
| Chain U(U(U(U(U(0, 5666))))) | **34.0%** | **22.7%** | **17.5%** | **13.8%** | **12.0%** |
| **(NON-Logarithmic)** US County Area | 38.5% | 14.5% | 17.2% | 13.5% | 16.3% |
| **(NON-Logarithmic)** Payroll Data | 45.6% | 6.6% | 8.5% | 18.3% | 21.0% |
| **(NON-Logarithmic)** Normal( 177, 40) | 39.2% | 2.5% | 8.4% | 20.4% | 29.5% |
| **(NON-Logarithmic)** Uniform(5, 78000) | 21.0% | 21.3% | 20.6% | 21.4% | 15.7% |
| **(NON-Logarithmic)** k/x over (1, 10) | 31.0% | 22.5% | 17.9% | 15.0% | 13.5% |

**FIGURE 5** 5-Bin Scheme - Arbitrarily Varying F (Start=0 Width=0.007)

To emphasize that only (digitally) logarithmic data types obey bin laws (patterns), two additional non-logarithmic data sets are added at the bottom, US County Area and Payroll data, serving as a contrast. These two data sets are well-known as being non-logarithmic (suffering from low order of magnitude). In addition, three non-logarithmic distributions are added, the Normal, the Uniform, and k/x defined over the short range (1, 10) which is non-logarithmic unless applied to digital BL base 10 or similarly constructed bin scheme. Results here demonstrate that the property of being either logarithmic or non-logarithmic is universal, number system invariant and quantitative in essence.

Even more surprising, arbitrary fractional values for the inflation factors vector Fi yield consistent results. Figure 6 depicts the results of a 6-bin scheme, starting at 0, with an initial width of 0.037, and utilizing the arbitrarily chosen set of Fi inflation fractional factors: {2.37, 3.08, 1.55, 4.17, 1.18, 2.35, 1.82, 5.07, 3.39, 2.04, 4.82, 7.07, 2.33, 6.67, 3.01, 1.67, 2.97, 3.33, 6.08, 2.25}. The near steady proportions here strongly suggests that there is no need whatsoever to fix Fi as integers in order to observe a common rate of the fall in histograms for all logarithmic data types! It should be emphasized that by now, especially with this last fractional Fi scheme, we have strayed quite far away from any connection to number system and digits, yet we are still able to obtain consistent and reliable quantitative laws.

| Data Set | Bin A | Bin B | Bin C | Bin D | Bin E | Bin F |
|---|---|---|---|---|---|---|
| Earthquakes | **27.7%** | **19.5%** | **16.3%** | **13.5%** | **12.3%** | **10.7%** |
| USA Population | **28.3%** | **20.2%** | **16.6%** | **14.0%** | **11.1%** | **9.8%** |
| Symmetrical Triangular | **29.0%** | **20.0%** | **15.8%** | **13.2%** | **11.7%** | **10.3%** |
| k/x over (1, 1000000) | **26.7%** | **20.8%** | **16.6%** | **13.8%** | **12.4%** | **9.8%** |
| Exponential Growth | **27.1%** | **19.8%** | **16.0%** | **14.1%** | **12.4%** | **10.7%** |
| Lognormal, L=5, S=1 | **27.6%** | **20.4%** | **16.3%** | **13.3%** | **11.7%** | **10.7%** |
| Lognormal, L=9.3, S=1.7 | **31.5%** | **20.5%** | **15.8%** | **12.3%** | **10.7%** | **9.1%** |
| Varied Data - Hill's Model | **28.2%** | **19.3%** | **16.7%** | **13.7%** | **12.5%** | **9.6%** |
| Chain 5 Uniforms | **26.0%** | **20.4%** | **16.1%** | **14.3%** | **12.4%** | **10.8%** |
| **(NON-Logarithmic)** US County Area | 26.6% | 25.0% | 16.5% | 15.0% | 9.3% | 7.6% |
| **(NON-Logarithmic)** Payroll Data | 31.3% | 10.8% | 14.4% | 15.1% | 15.3% | 13.1% |
| **(NON-Logarithmic)** Normal(177, 40) | 23.9% | 42.9% | 24.8% | 5.2% | 1.5% | 1.8% |
| **(NON-Logarithmic)** Uniform(5, 78000) | 19.3% | 16.0% | 16.3% | 15.7% | 16.5% | 16.2% |
| **(NON-Logarithmic)** k/x over (1, 10) | 21.1% | 26.7% | 21.6% | 11.3% | 10.0% | 9.3% |

**FIGURE 6** 6-Bin Scheme - Fractional Arbitrarily Varying F (Start=0 Width=0.037)

Interestingly, we can put even those bin schemes of arbitrarily varying inflation factors on (some) 'theoretical' basis as a confirmation of sorts by considering the average value of those varying Fi inflation values, thus 'enabling' ourselves to insert this singular $F_{AVG}$ value within the algebraic series of infinitely expanding k/x. $F_{AVG}$ value for the integral arbitrary factors of Figure 5 is 4.55, while $F_{AVG}$ value for the fractional arbitrary factors of Figure 6 is 3.36. Performing exactly this comparison (empirical to 'theoretical') for these two bin schemes we obtained the following results:

Arbitrarily chosen integral Fi values (Figure 5):
Average of empirical bin results of 9 data sets: {35.6%, 22.7%, 17.0%, 13.3%, 11.4%}
'Theoretical' k/x infinitely expd. D=5, F= 4.55: {35.5%, 22.9%, 17.0%, 13.5%, 11.2%}

Arbitrarily chosen fractional Fi values (Figure 6):
Avg empirical bin results of 9 data sets: {28.0%, 20.3%, 16.2%, 13.5%, 11.8%, 10.1%}
'Theoretical' k/x inf. expd. D=6, F= 3.36: {27.5%, 20.5%, 16.4%, 13.6%, 11.7%, 10.2%}

In both bin schemes, the 'theoretical' and empirical results are quite close, and in spite of the fact that $F_{AVG}$ value is just a 'simple-minded' average of the varying inflation factors. This certainly lends varying F factors bin schemes more 'theoretical respectability'. Surely there exists no fractional base within a number system, nor any such concept as seven and a half bin scheme, both must be whole numbers, but fractional inflation factors are totally valid conceptual possibilities, even without the above agreement between empirical and 'theoretical' results. One must bear in mind that in the derivation of the series of algebraic expressions for the infinitely expanded k/x distribution no assumption was made restricting F to an integral value, although there is no mathematical basis at all for employing and extrapolating the result of the infinite x/k model with steady F to a model with varying Fi (in spite of the apparent 'empirical' near correspondence in numerical results above).

## [13] Higher Orders Digits in BL Interpreted as Particular Bin Schemes

2nd order digit distribution in base 10 can be viewed as the arrangement of equally-spaced 10 inner bins within each outer 1st order bin. These inner bins come with a halting/vacillating process of expansion. If presented simply in terms of a bin scheme with varying inflation factors Fi, and focusing only on the x-axis part over 1 for brevity, the vector of Fi in a 10-bin scheme tailor-made for 2nd order digits is as follow:
Fi ={1,1,1,1,1,1,1,1,**10**,1,1,1,1,1,1,1,1,**10**,1,1,1,1,1,1,1,1,**10**,1,1,1,1,1,1,1,1,**10**,1,1,1, etc.}.
Simple-mindedly extrapolating the mathematical results of k/x infinitely expanded with the singular insertion of $F_{AVG}$ value in such bin system interpretation of the second order with $F_{AVG}$ value of (1+1+1+1+1+1+1+1+10)/9 = 18/9 = 2, and D=10, we obtain:

Infinite k/x Bin Scheme D=10 F=2: {13.8, 12.6, 11.5, 10.7, 10.0, 9.3, 8.7, 8.2, 7.8, 7.4}
Digital BL 2nd Order Distribution: {12.0, 11.4, 10.9, 10.4, 10.0, 9.7, 9.3, 9.0, 8.8, 8.5}

The two results are a bit close, but not identical. If instead of 2 for F the value 1.49106 is employed, the two results are extremely close to each other, yet such empirical explorations are by no means an appropriate substitution for rigorous mathematics!

Infinite k/x Bin D=10 F=1.49106: {12.0, 11.5, 10.9, 10.5, 10.1, 9.7, 9.3, 9.0, 8.7, 8.4}
Digital BL 2nd Order Distribution: {12.0, 11.4, 10.9, 10.4, 10.0, 9.7, 9.3, 9.0, 8.8, 8.5}

To match digital 2nd order distribution in BL with our bin theory, it is necessary to skip the theoretical infinitely expanded k/x model altogether (as it contains no mathematical results for varying Fi factors) and simply perform an empirical bin scheme with D=10 and varying F factors Fi = {1,1,1,1,1,1,1,1,**10**,1,1,1,1,1,1,1,1,**10**,1,1,1, etc.}. The results of 3 such schemes with a variety of initial widths and starting points on the 2009 US population data sets are as follow:

Starting at 0.5, initial width 0.07: {12.6, 11.7, 10.3, 10.5, 10.3, 9.4, 9.1, 9.2, 8.9, 8.0}
Starting at 0, initial width 0.30: {12.3, 11.4, 10.5, 10.1, 10.2, 9.5, 9.7, 9.0, 9.2, 8.2}
Starting at 0.63, width 0.00045: {11.9, 11.0, 11.4, 10.2, 10.3, 9.7, 9.6, 8.6, 9.2, 8.1}
Digital BL 2nd Order Distribution: {**12.0, 11.4, 10.9, 10.4, 10.0, 9.7, 9.3, 9.0, 8.8, 8.5**}

Applying such an interpretation, two explanations can be made of why 2nd order digit distribution in BL is more equal and less skewed in comparison with 1st digit distribution: (I) The widths of the inner bins of the 2nd order are much smaller and more refine than those of the outer bins of the 1st order (on any given segment of the x-axis), hence fall in density is registered as less severe. (II) Inner bins are not constantly expanding, while outer bins are constantly expanding in each cycle. The ramification of such bin vista for higher orders digits in general is quite profound! It implies that any well-defined bin structure which alternates between flatness and expansion gives birth to a bin law such that for a given (fixed) number of bins, resultant bin skew-ness depends on degree of expansion, namely a one-to-one relationship between skew-ness and degree of expansion. Flat bin schemes yield bin equality; fast-expanding bin schemes yield extreme bin inequality (super skew-ness); while schemes with mild flatness having some expansion yield intermediate results.

This crucial and fascinating result relates to one of the most central issue in BL and bin systems in general, namely that data falls on the x-axis in a multiplicative-like manner, just as seen in Multiplicative CLT models. That there is acceleration in how the spread of the data gets diluted towards the right side of high values, which can only be seen and detected via a bin prism of just as rapidly expanding multiplicative inflation factors! As seen here, flat bin schemes (F=1) are unable to convey any numerical measure of fall, but rather always yield bin-equality by default whenever width is sufficiently refine, for all data types, logarithmic or otherwise, and as such do not constituted a law.

As one demonstration, several 7-bin schemes are performed on 2009 US Population Centers Census data from the origin with initial width 0.0039, each with a different F expansion factor. Results are shown in Figure 7.

| F Inflation | Bin A | Bin B | Bin C | Bin D | Bin E | Bin F | Bin G |
|---|---|---|---|---|---|---|---|
| 1 | 14.6% | 14.2% | 14.2% | 13.9% | 14.2% | 14.2% | 14.6% |
| 2 | 19.2% | 17.6% | 14.6% | 13.9% | 12.1% | 11.8% | 10.7% |
| 3 | 23.2% | 17.7% | 14.9% | 12.7% | 11.7% | 10.5% | 9.4% |
| 4 | 25.5% | 19.0% | 15.1% | 12.2% | 10.4% | 9.2% | 8.5% |
| 5 | 28.2% | 19.2% | 14.5% | 12.0% | 9.7% | 9.0% | 7.5% |
| 6 | 30.0% | 19.5% | 14.5% | 11.4% | 9.2% | 8.0% | 7.3% |
| 7 | 32.2% | 19.2% | 14.4% | 11.2% | 8.7% | 7.7% | 6.7% |
| 8 | 33.4% | 19.2% | 13.6% | 10.8% | 9.0% | 7.5% | 6.6% |
| 9 | 34.5% | 19.4% | 14.1% | 10.5% | 8.5% | 6.9% | 6.0% |
| 10 | 35.6% | 19.8% | 13.2% | 9.8% | 8.5% | 7.1% | 6.0% |
| 11 | 37.7% | 18.9% | 13.1% | 10.1% | 8.1% | 6.3% | 5.8% |
| 12 | 38.1% | 19.1% | 13.1% | 9.5% | 7.8% | 6.6% | 5.8% |

**FIGURE 7** Bin-skew-ness Increases in Direct Proportion to Inflation Factor F – US Population Centers Data

This decisive demonstration can **not** be performed in the context of digital BL where variations in F and D go hand in hand, since one is not allowed to vary independently of the other (they must relate to each other as in F = D + 1). It is only in the context of the bin theory developed here that we can isolate F and exclusively vary it alone; able to observe this relationship and dependency; and thus demonstrate the principle.

As a general conclusion, Base and Order in digital BL are mere variations on bin scheme structure. That is, they are simply different parameters within bin schemes. Digital **Base** in the context of bin theory is the number of bins plus 1 and it is expressed as (D+1), while digital **Order** in the context of bin theory is simply a particular combinations of varying Fi values inflation factors.

The dependency of flat bin schemes (F = 1) on the value of the width shall be demonstrated by examining 2009 Census data on US Population centers. The table in Figure 8 depicts a wide variety of results and skew-ness conditions of 4-bin schemes depending on the value of the width. The starting point clearly does not play any role here unless it is placed very far from the origin so that a significant portion of the data escapes the scrutiny of the applied bin scheme. For width size that is too small compared with the outlay of this population data, bin sizes are simply too refined to show any meaningful differentiation in data concentration, therefore discrete values fall upon the bins in a totally random manner resulting in bin equality. For width size that is too big compared to the outlay of this population data, bin sizes are simply too crude, capturing significant portions of the data within the first bins of the first cycles distorting their message of relative concentration. The huge width size of 10,000,000 captures the entire data within the first bin of the first cycle (allocating it 100%), and renders the whole scheme meaningless.

Finally, for the infinitely expanded flat bin schemes constructed for the k/x distribution (F= 1), the theoretical result should yield bin-equality (that is: 1/D per bin) since the curve in the limit 'towards infinity' on the far right 'appears flat' and thus intuitively distributes equal proportions to all bins.  When 1 is substituted for F in the expression of the infinite sequence, convergence is not apparent at all in computer simulations/calculations for lack of 'infinite' time and 'limitless' supply of memory chip units.

| Width | Starting | Bin A | Bin B | Bin C | Bin D |
|---|---|---|---|---|---|
| 1 | 1 | 25.9% | 25.2% | 24.2% | 24.7% |
| 3 | 1 | 25.0% | 25.6% | 24.5% | 24.9% |
| 20 | 1 | 25.0% | 25.3% | 25.0% | 24.7% |
| 50 | 1 | 25.0% | 26.0% | 25.4% | 23.7% |
| 50 | 2 | 25.1% | 26.1% | 25.3% | 23.6% |
| 50 | 3 | 25.2% | 25.8% | 25.4% | 23.6% |
| 100 | 1 | 26.2% | 26.9% | 24.8% | 22.1% |
| 135 | 1 | 27.7% | 27.8% | 23.6% | 20.8% |
| 170 | 1 | 29.2% | 28.1% | 22.8% | 20.0% |
| 250 | 0 | 33.1% | 27.6% | 21.6% | 17.7% |
| 250 | 1 | 33.2% | 27.5% | 21.6% | 17.7% |
| 400 | 1 | 39.5% | 26.6% | 18.8% | 15.1% |
| 700 | 1 | 49.3% | 23.4% | 15.3% | 12.0% |
| 700 | 13 | 49.5% | 23.3% | 15.2% | 11.9% |
| 700 | 19 | 49.7% | 23.1% | 15.2% | 12.0% |
| 1200 | 1 | 58.8% | 19.9% | 12.5% | 8.9% |
| 2500 | 1 | 70.6% | 14.8% | 8.7% | 6.0% |
| 10 Million | 1 | 100.0% | 0.0% | 0.0% | 0.0% |

**FIGURE 8** Results of Flat Bins (F = 1) Depends on Width - US Population Centers

## [14] Bin Development Pattern

A bin developmental pattern along bin cycles is certainly expected to emerge upon careful examination (in all random data types), and regardless of expansion style. Within the first cycle of D bins, approximate bin equality may prevail, and even some bin-reversal of fortunes may occur. Around the central cycles where most of the data resides one should find bin proportions that closely match the relevant bin-law given parameters D and F. Finally, around the far right region, more extreme bin-inequality (super skew-ness) should prevail. A related article by this author on BL demonstrates this principle (termed 'digital development pattern') which could only be applied for the particular partition of the x-axis along sub-intervals standing between integral powers of ten, where on the socialist left region an approximate digital equally and harmony prevail; around the center mini digit distributions closely mimic $LOG_{10}(1+1/d)$; and on the extreme far right there exist severe and harsh digital inequality in favor of low digits (i.e. skew-ness over and above the logarithmic distribution). Such digital development was clearly demonstrated in the generic cases of the Lognormal and the Exponential distributions, but it is totally absent in the case of the k/x distribution (defined over a very long range, preferably with an integral exponent difference between its edges). The case of such k/x distribution shows a steady and consistent $LOG_{10}(1+1/d)$ behavior throughout its entire range.

Certainly all this variation in skew-ness **within** a single bin scheme (i.e. development) is in perfect harmony and nicely consistent with the variation seen in skew-ness **between** different schemes having a variety of F inflation factors as seen Figure 7. This is so since even a fast expanding scheme with high F factor achieves its widest bin sizes on the right after numerous application of such F factor; while in the beginning on the left, it must start with a very narrow and refined width, and it is yet to apply sufficient number of such F multiplications. It should be noted that bin development pattern exists no matter how large an inflation factor F is chosen relative to D bin number. However, high inflation factors F relative to D are associated with faster and more dramatic bin development. Low inflation factors F relative to D are associated with gradual and less dramatic bin development. Data congregation however always occurs around the central bins, and very little data is found on the left-most or right-most bins, regardless of bin style.

Moreover, all this is also nicely consistent with the intuition that [given a particular falling density curve] too thick and crude a net of bins cast upon the x-axis is associated with over-skew-ness as it records **macro** fall, while a refined and thin net is associated with less skew-ness as it records **micro** fall.

But in the generic case of (logarithmic) k/x defined over a very large range, bin development depends on expansion style. For a number-system-like expansion style with F = D + 1, bin proportions are steady and consistent throughout the entire bin-cycle structure, with each bin cycle containing the same amount of overall data, just as was seen for our number system with base 10 in terms of digits and BL. For a fast expanding bin style where inflation factor F > D + 1, bin proportions show development towards skewer distributions, reaching a certain skew-ness level very rapidly, and then maintaining and continuing it thereafter in all subsequent cycles. Also data portion within each bin cycle increases, but rapidly reaches a certain level, maintaining it about constant thereafter in all subsequent cycles. For a slow expanding bin style where inflation factor F < D + 1, bin proportions show an inverse development towards less skewed and more even distributions, rapidly reaching a certain (milder) skew-ness level, and then maintaining and continuing it thereafter in all subsequent cycles. Also data portion within each bin cycle decreases, but rapidly reaches a certain level, then maintaining it about constant thereafter in all subsequent cycles.

All this is consistent with the theoretical expressions of the series of infinitely expanding bin cycles on k/x derived earlier. In fact is can be mathematically deduced from it. In the theoretical k/x case where F = D +1, once-expanding and twice-expanding schemes were already shown to reduced to the same proportion of LOG(1+1/d)/LOG(D+1) corresponding to non-expanding bin scheme. This reduction was assumed to hold for any N-expanding scheme, thus implying a constant bin proportions throughout all expanding cycles. This is so since a single expansion of k/x (say) does not alter bin proportions of the original non-expanding x-axis range, where value of k is admittedly reduced due to the expansion of the defined range [entire area must always sum to 1.] This reduction in the value of k affects height of k/x over original non-expanding range, but not bin proportions there since area of each bin is reduces by an identical factor. Hence, if overall bin proportions stay unaffected with each expansion, it follows that with each expansion an extra set of identical bin proportions must be added (i.e. no bin development pattern). This is analogous to numbers being added one at a time where the average keeps its value as a constant, implying that all the numbers being added are identical (having the value of the average). Similar mathematical arguments can be made regarding F > D + 1 scenario, showing that bin proportions there must be changing and becoming skewer with each expansion. Same argument can be made showing that milder (less skewed) proportions are gotten with each expansion in the F < D + 1 scenario.

## [15] The General Scale Invariance Principle

The scale invariance principle, originally stated in the context of digits and BL, is a much more general principle, extendable to all bin schemes, namely that the vector of proportions of **any** given logarithmic data set X examined via **any** given bin scheme is the same as the vector of proportions of K*X examined via the same bin scheme for K being **any** positive real number.

Empirical examination of bin results for a variety of re-scaling transformations of several physical data set confirms the general status of the scale invariance principle with an steady resultant bin proportions no matter the scale in use. Empirical examination of bin results for a variety of multiplicative transformations of realizations from k/x distribution defined over (14, 231316943) also confirms the general status of the scale invariance principle in this theoretical case.

Let us mathematically prove the general principle of scale invariance in the two most crucial and relevant distributions within the context of BL; namely for k/x distribution as well as for the Lognormal distribution. In both cases the Distribution Function Technique shall be applied using somewhat similar notations as in Freund's book "Mathematical Statistics" - sixth edition chapter 7.3. We shall assume that the data set X (actualizations of real values from the distribution) is being transformed (re-scaled) by the factor K into Y = K*X.

In the distribution case pdf(X) = k/x over (a, b), the monotonic transformation
Y(X) = K*X implies that X(Y) = Y/K, and dx/dy = 1/K, hence:

**pdf(Y) = pdf( X(Y) ) * | dx/dy |**
pdf(Y) = K/[Y/K] *| 1/K | = K/Y over (aK, bK),

Namely, multiplying k/x-like data set by any constant K yields the same density form, although the range is different. Moreover, the change in range should not adversely affect logarithmic behavior whatsoever, and to demonstrate that we shall consider two contexts:

(I) A short (a, b) range in the context of digits and BL with base B, and where a and b spans exactly an integral exponent difference of the base, namely that $LOG_B(b) - LOG_B(a) = INTEGER$.
(II) The generic definition of logarithmic-ness for k/x defined over a truly huge range.

In the first context, the newly defined range is still of an integral exponent difference, since $LOG_B(Kb) - LOG_B(Ka) = LOG_B(K) + LOG_B(b) - LOG_B(K) - LOG_B(a) = LOG_B(b) - LOG(a)$, which is also an INTEGER as was assumed for the original data set.

In the second context, the newly defined range is still as huge if K=1; much longer if K>1; but even when K<1 the true measure of 'length' is log-wise which is unchanged.

Let us turn our attention now to the case of the Lognormal distribution, namely to $pdf(X) = [1/(X\sigma\sqrt{(2\pi)})]*e$ to $(-1/2)*[(\ln(X) - \mu)/\sigma]^2$. The monotonic transformation $Y(X) = K*X$ implies that $X(Y) = Y/K$, and $dx/dy = 1/K$, hence:

**pdf(Y) = pdf( X(Y) ) * | dx/dy |**
$pdf(Y) = [1/((Y/K)\sqrt{(2\pi)})]*e$ to $(-1/2)*[(\ln(Y/K) - \mu)/\sigma]^2 * | 1/K |$
The values of K on the left and on the right cancel out, and we are left with:
$pdf(Y) = [1/((Y)\sigma\sqrt{(2\pi)})]*e$ to $(-1/2)*[(\ln(Y) - \ln(K) - \mu)/\sigma]^2$
This can be rewritten as:
$pdf(Y) = [1/(Y\sigma\sqrt{(2\pi)})]*e$ to $(-1/2)*[(\ln(Y) - \mathbf{(\ln(K) + \mu)})/\sigma]^2$
Namely another Lognormal distribution of the same shape parameter σ, but of a different location parameter ln(K) + μ.

Since logarithmic behavior of the lognormal is determined solely by the value of the shape parameter (being logarithmic whenever shape is roughly larger than 1), the fact that shape parameter is totally unaffected under re-scaling confirms the scale invariance principle for this crucial and very relevant statistical distribution.

As discussed in earlier chapters, logarithmic-ness is an absolute and universal property of any data set or distribution, irrespective of number system in use, base, or bin scheme chosen. Consequently, guaranteeing logarithmic behavior (in the context of digits and BL) for any re-scaled Lognormal and k/x distribution implies guaranteeing their compliance with any bin scheme chosen, namely the general scale invariance principle for bin systems as well (in those two cases at least).

## [16] Closed Form Expression for the Limit of the Infinite Sequence

Let us calculate an analytical (closed form) expression for the limit of the infinite sequence of k/x bin scheme model, enabling us to succinctly express the general law of relative quantities. This is calculated first for the F > 1 case. The assumption F ≠ 1 is solely an initial restriction enabling progress in the reduction, and it shall be relaxed later.

The 4th term of the sequence expressed earlier, denoted as $S_4$ is:

$$S_4 = \frac{\ln\left(\frac{[1+(d)]}{[1+(d-1)]}\right) + \ln\left(\frac{[1+D+(d)F]}{[1+D+(d-1)F]}\right) + \ln\left(\frac{[1+D+DF+(d)F^2]}{[1+D+DF+(d-1)F^2]}\right) + \ln\left(\frac{[1+D+DF+DF^2+(d)F^3]}{[1+D+DF+DF^2+(d-1)F^3]}\right)}{\ln(1 + D + DF + DF^2 + DF^3)}$$

Employing the **finite** geometric formula for the terms involving F, namely:
$1 + X + X^2 + X^3 + \ldots + X^N = (X^{N+1} - 1) / (X - 1)$ ,
the nth term in the sequence is then:

$$S_N = \frac{\displaystyle\sum_{j=0}^{j=N-1} \ln\left(\frac{1 + \frac{D(F^j - 1)}{(F-1)} + (d)F^j}{1 + \frac{D(F^j - 1)}{(F-1)} + (d-1)F^j}\right)}{\ln\left(1 + \frac{D(F^N - 1)}{(F-1)}\right)}$$

Pulling together all the coefficients of $F^{POWER}$, we get:

$$S_N = \frac{\displaystyle\sum_{j=0}^{j=N-1} \ln\left(\frac{1+\left(\frac{D+(d)(F-1)}{(F-1)}\right)F^j - \frac{D}{F-1}}{1+\left(\frac{D+(d-1)(F-1)}{(F-1)}\right)F^j - \frac{D}{F-1}}\right)}{ln\left(1+\left(\frac{D}{(F-1)}\right)F^N - \frac{D}{F-1}\right)}$$

In order to obtain a more compact expression, let us define:

$$A = \frac{D + d(F-1)}{F-1}$$

$$B = \frac{D + (d-1)(F-1)}{F-1}$$

$$C = \frac{D}{F-1}$$

$$E = 1 - \frac{D}{F-1}$$

$$S_N = \frac{\displaystyle\sum_{j=0}^{j=N-1} \ln\left(\frac{AF^j + E}{BF^j + E}\right)}{ln(CF^N + E)}$$

Since in our context $F \geq 1$, namely either 1 as in flat bin schemes, or larger than 1 as in normal expanding bin schemes yielding quantitative laws, there is no hope of obtaining any obvious convergence in terms such as $F^j$ or $F^N$, hence we define $f = 1 / F$, creating a quantity f such that $0 < f \leq 1$ holds, and which may hopefully let terms such as $f^j$ or $f^N$ converge.

$$S_N = \frac{\displaystyle\sum_{j=0}^{j=N-1} \ln\left(\frac{A\frac{1}{f^j}+E}{B\frac{1}{f^j}+\ E}\right)}{ln\left(\ C\frac{1}{f^N}+E\ \right)}$$

$$S_N = \frac{\displaystyle\sum_{j=0}^{j=N-1} \ln\left(\frac{A\frac{1}{f^j}+E}{B\frac{1}{f^j}+\ E} *\ \frac{f^j}{f^j}\right)}{ln\left(\ C*\frac{1}{f^N}+E\ \ *\ \frac{C}{C}*\ \frac{f^N}{f^N}\right)}$$

$$S_N = \frac{\displaystyle\sum_{j=0}^{j=N-1} \ln\left(\frac{A+E*f^j}{B+\ E*f^j}\right)}{ln\left(\ C*\frac{1}{f^N}*\left(\ 1+E\ \ *\ \frac{f^N}{C}\right)\right)}$$

$$S_N = \frac{\displaystyle\sum_{j=0}^{j=N-1} \ln\left(\frac{\left(1+\frac{E}{A}*f^j\right)A}{\left(1+\ \frac{E}{B}*f^j\right)B}\right)}{ln\left(\ C*F^N*\left(\ 1+\ \frac{E}{C}f^N\right)\right)}$$

$$S_N = \frac{N * ln\left(\frac{A}{B}\right) + ln\left(\prod_{j=0}^{N-1} \frac{\left(1 + \frac{E}{A} * f^j\right)}{\left(1 + \frac{E}{B} * f^j\right)}\right)}{N * ln(F) + ln(C) + ln\left(1 + \frac{E}{C} f^N\right)}$$

It is only at this late stage that we let N go to infinity!
**For F > 1 in the normal case of expanding bin scheme**, 0 < f < 1, therefore

$$\prod_{j=0}^{\infty} \frac{\left(1 + \frac{E}{A} * f^j\right)}{\left(1 + \frac{E}{B} * f^j\right)}$$

is a convergent infinite product since $\sum_{j=0}^{\infty} f^j$ is converging. The term $ln\left(1 + \frac{E}{C} f^N\right)$ is zero as $N \to \infty$. The term ln(C) is ln($\frac{D}{F-1}$), and it is also finite as $N \to \infty$, finally:

$N \to \infty \quad S_N = \frac{\ln(\frac{A}{B})}{\ln(F)}$ and using the definition of A and B above, the general relative quantities law is then $\frac{\ln\left(\frac{\left(\frac{D+d(F-1)}{F-1}\right)}{\left(\frac{D+(d-1)(F-1)}{F-1}\right)}\right)}{\ln(F)}$, which is further reduced by canceling out the two (F – 1) terms in the numerator to arrive at:

$$\text{The General Law:}\quad \frac{\ln\left(\dfrac{D + d(F-1)}{D + (d-1)(F-1)}\right)}{\ln(F)}$$

To verify that digital BL is simply a special case and a consequence of the general law of relative quantities when bin schemes are constructed under the constraint F = D + 1, the term F is then substituted by D + 1 everywhere in expression of the general law:

$$\text{GL} = \frac{\ln\left(\dfrac{D+d(D+1-1)}{D + (d-1)(D+1-1)}\right)}{\ln(D+1)} = \frac{\ln\left(\dfrac{D+d(D)}{D + (d-1)(D)}\right)}{\ln(D+1)} = \frac{\ln\left(\dfrac{1+d}{1+(d-1)}\right)}{\ln(D+1)}$$

$$= \frac{\ln\left(1+\frac{1}{d}\right)}{\ln(D+1)} = \frac{\ln\left(1+\frac{1}{d}\right)}{\ln(BASE)} = \frac{\text{LOG}\left(1+\frac{1}{d}\right)}{\text{LOG}(BASE)} = \text{BL}$$

For a number system with base 10, we get:

$$\text{BL} = \frac{\text{LOG}\left(1+\frac{1}{d}\right)}{\text{LOG}(10)} = \frac{\text{LOG}\left(1+\frac{1}{d}\right)}{1} = \text{LOG}_{10}\left(1 + \frac{1}{d}\right)$$

## [17] Closed-Form Expression for the Limit in the Case F = 1

In the case of flat bin schemes with **F = 1**, the 4th term of the sequence denoted as $S_4$ is:

$$S_4 = \frac{\ln\left(\frac{[1+(d)]}{[1+(d-1)]}\right) + \ln\left(\frac{[1+D+(d)]}{[1+D+(d-1)]}\right) + \ln\left(\frac{[1+D+D+(d)]}{[1+D+D+(d-1)]}\right) + \ln\left(\frac{[1+D+D+D+(d)]}{[1+D+D+D+(d-1)]}\right)}{\ln(1+D+D+D+D)}$$

$$S_N = \frac{\sum_{j=0}^{j=N-1} \ln\left(\frac{1+D*j+(d)}{1+D*j+(d-1)}\right)}{\ln(1+N*D)}$$

$$S_N = \frac{\sum_{j=0}^{j=N-1} \ln\left(\frac{1+D*j+d}{D*j+d}\right)}{\ln\left(N*\left(1/N+D\right)\right)}$$

As N → ∞, the <u>numerator</u> (NU) can be roughly evaluated by using the Integral Test Approximation.

$$\text{NU} = \int_0^N \; ln\left(\frac{1+Dx+d}{Dx+d}\right) dx =$$

$$= \int_0^N [\; ln(\,1 + Dx + d) - ln(\,Dx + d)]\; dx$$

$$= \int_0^N ln(\,1 + Dx + d)\, dx \;-\; \int_0^N ln(\,Dx + d)\, dx$$

Let $u = (1 + Dx + d)$ , implying that $du/dx = D$

Let $z = (Dx + d)$ , implying that $dz/dx = D$

$$\text{NU} = \int_{1+d}^{1+DN+d} ln(\,u)\frac{1}{D}du \;-\; \int_{d}^{ND+d} ln(\,z)\frac{1}{D}dz$$

Using $\int \ln(x)\, dx = x * ln(x) - x + C$ , we evaluate the numerator as:

-----------------------------------------------------------------------------------------

$$+\frac{1}{D}(1 + DN + d)\ln(1 + DN + d) - \frac{1}{D}(1 + DN + d)$$

$$-\frac{1}{D}(1 + d)\ln(1 + d) + \frac{1}{D}(1 + d)$$

$$-\frac{1}{D}(DN + d)\ln(DN + d) + \frac{1}{D}(DN + d)$$

$$+\frac{1}{D}(d)\ln(d) - \frac{1}{D}(d)$$

-----------------------------------------------------------------------------------------

The terms not involving N at all (in gray color) are negligible in the limit as N goes to infinity and can be omitted, hence we get:

$$+\frac{1}{D}(1 + DN + d)\ln(1 + DN + d) - \frac{1}{D}(DN)$$

$$-\frac{1}{D}(DN + d)\ln(DN + d) + \frac{1}{D}(DN)$$

-----------------------------------------------------------------------------------------

$$+\frac{1}{D}(1 + DN + d)\ln(1 + DN + d) \;-\; \frac{1}{D}(DN + d)\ln(DN + d)$$

-----------------------------------------------------------------------------------------

-----------------------------------------------------------------------------------------

$\frac{1}{D} * [(1)\ln(1 + DN + d) + (DN + d)\ln(1 + DN + d) - (DN + d)\ln(DN + d)]$

-----------------------------------------------------------------------------------------

$$\frac{1}{D} * [\ \ln(DN)\ +\ (DN + d) * \ln\left( {(1 + DN + d)}\big/{(DN + d)} \right)\ ]$$

-----------------------------------------------------------------------------------------

$$\frac{1}{D} * [\ \ln(D) + \ln(N) + (DN + d) * \ln\left( 1 + {1}\big/{(DN + d)} \right)\ ]$$

-----------------------------------------------------------------------------------------

$$\mathbf{NU} \ = \ \frac{1}{D} * [\ \ln(N) + \ln\left( \left( 1 + {1}\big/{(DN + d)} \right)^{(DN+d)} \right)\ ]$$

-----------------------------------------------------------------------------------------

The inner expression inside the natural logarithm on the right is simply Euler's number $\mathbf{e}$ (the exponential constant) defined as $\mathbf{lim}_{n \to \infty} \left(\mathbf{1} + \frac{\mathbf{1}}{\boldsymbol{n}}\right)^{\boldsymbol{n}}$ which converges to the finite 2.71828 value. Hence the whole ratio of the numerator divided by the denominator in the limit as N approaches infinity is:

$$\mathbf{S}_{\mathbf{N} \to \infty} = \frac{\frac{1}{D}[\ \ln(N) + \ln(e)]}{\ln(N) + \ \ln\left( {1}\big/{N} + \ D \right)}$$

$$\mathbf{S}_{\mathbf{N} \to \infty} = \frac{\left(\frac{1}{D}\right)\ln(\mathrm{N})}{\ln\ (\mathrm{N})} = \frac{\mathbf{1}}{\mathbf{D}}$$

Namely bin equality for flat bin schemes with F = 1, regardless what value D takes.
The mathematician George Andrews [who also assisted in the reduction in the F > 1 case] suggests an alternative and concise proof in the F = 1 case, derived directly from the general result of the proof in the F > 1 case via L'Hopital's Rule, and in spite of the fact that F ≠ 1 was a necessary assumption in parts of the derivation there [such as A, B, C, and E being undefined when F is 1, and violating the geometric formula restriction on common ratio not assuming the value of 1]. Andrews' approach is to let N approach infinity first as in the general F ≠ 1 proof, and only then let F approach 1 in the expression of the General Law. Direct substitution by 1 for F in the General Law yields the indeterminate form:

$$\frac{\ln\left(\frac{D+d(1-1)}{D+(d-1)(1-1)}\right)}{\ln(1)} = \frac{\ln\left(\frac{D+0}{D+0}\right)}{\ln(1)} = \frac{\ln(1)}{\ln(1)} = \frac{0}{0}$$

$$\lim_{F\to 1} \frac{\ln\left(\frac{D+d(F-1)}{D+(d-1)(F-1)}\right)}{\ln(F)} =$$

$$\lim_{F\to 1} \frac{\ln(D+d(F-1)) - \ln(D+(d-1)(F-1))}{\ln(F)} =$$

Applying L'Hopital's Rule and differentiating with respect to F, we obtain:

$$\lim_{F\to 1} \frac{\frac{d}{(D+d(F-1))} - \frac{(d-1)}{(D+(d-1)(F-1))}}{1/F} =$$

Direct insertion of 1 for F yields:

$$\frac{\frac{d}{(D+d(1-1))} - \frac{(d-1)}{(D+(d-1)(1-1))}}{1/1} =$$

$$\frac{d}{(D+d(0))} - \frac{(d-1)}{(D+(d-1)(0))} = \frac{d}{D} - \frac{(d-1)}{D} = \frac{1}{D}$$

This latter proof is about the limit as F approaches 1, while the former proof is about F actually attaining the exact value of 1. Surely one cannot take for granted that the limit of f(x) as x approaches K equals f(K) since a discontinuity might be a distinct possibility. Therefore on the face of it, the two proofs which yield the same 1/D result and complement each other, are actually two distinct and necessary results. Yet, as Andrews further points out, both the numerator and the denominator are analytic around F = 1 and L'Hopital's rule guarantees that each had a first order zero at F = 1. Consequently their quotient is analytic around F = 1 and its value at F = 1 is precisely the limit that was calculated. Hence the resultant bin equality of 1/D in his proof applies not only to the limit as F approaches 1, but also exactly at that point where F actually equals 1. This renders Andrews' proof not only more straightforward, but also entirely sufficient. Experimentations with two alternative versions of proofs with the aid of L'Hopital's rule by letting F approach 1 first and then letting N approach infinity have ended in decisive failures.

## [18] The Universal Law of Relative Quantities

Our general model measuring relative quantities via bin schemes that expand steadily by way of a constant F factor [either $F \neq D + 1$ or $F = D + 1$], leads to the G.L.O.R.Q., which is totally divorced from digits and independent of any number system. Yet this result does not tell the whole story of relative quantities in sufficient generality, as seen in chapter [12] earlier, where arbitrarily varying and fractional factors for F are employed. Since it was (successfully) postulated that the generic pattern in how relative quantities are found in nature is such that the frequency of quantitative occurrences is inversely proportional to quantity, leading to the k/x distribution as the mathematically fitting density, any arbitrarily-constructed bin scheme should point to a new law by way of the evaluation of definite integrals constructed for the k/x distribution defined over an infinite range of (w, +∞) utilizing exactly such a bin system. For example, for a bin scheme with a vector of expansion Fi = {1 + G/1, 1 + G/2, 1 + G/3, 1 + G/4, and so forth}, G being any positive real number, the **universal** law of relative quantities should be obtained mathematically by evaluating definite integrals for such a bin scheme constructed for k/x defined over an infinite range, leading to a different infinite algebraic series than the one in chapter [8] - and hopefully to a closed-form expression of that law if a limit of that series can be found. The above postulate about the generic pattern of relative quantities is nearly certain to be successfully applied to any bin scheme, and this conjecture is strongly supported by empirical confirmation of bin-behavior of k/x for varying and arbitrary fractional factors as in Figures 5 and 6, where k/x defined over the relatively large range of (1, 1,000,000) [albeit not infinite] corresponds nicely to the bin proportions of all the other logarithmic data sets. Several other empirical experimentations with odd and arbitrary bin schemes also yielded consistent results (i.e. bin laws) where k/x defined over (3, 734,336,442) also participated nicely in the group 'logarithmic' behavior.
Note: a large data set of 20,000 values was realized from k/x defined over (3, 734,336,442) by way of computer simulations. In general, for any k/x defined over $(10^A, 10^B)$, realizations of values are computer-simulated using the expression:

$$10^{A + \text{Uniform}(0, 1)*(B - A)}$$

## [19] Related Considerations

Single-issue physical data sets such as time between earthquakes, quasar rotation rates, population count, etc. pose a dilemma in the field of BL since they cannot be explained via the standard mathematical model of distribution of all distributions. It would be very hard or perhaps impossible to argue that measured time between earthquakes for example is some sort of a mixture of numerous distributions. Let us conjecture two alternative explanations here: (I) Multiplicative CLT leading to the Lognormal as the underling distribution, assuming that the physical process can be model approximately at least on some repeated multiplicative process. (II) Chain of distribution, assuming that the physical process can be modeled as having one of its parameter dependent on the random variable of another physical process, namely the interconnectedness and causality in life.

# Acknowledgment:

The author would like to thank the distinguished mathematician **George Andrews** of Pennsylvania State University for his assistance in the reduction of the infinite series for the bin systems constructed for k/x distribution, leading to a closed-form expression of The General Law of Relative Quantities.

Email: akossovs@yahoo.com
March 2014